\documentclass{article}

\usepackage[nonatbib,preprint]{nips_2018}

\usepackage[utf8]{inputenc} 
\usepackage[T1]{fontenc}    
\usepackage{hyperref}       
\usepackage{url}            
\usepackage{booktabs}       
\usepackage{amsfonts}       
\usepackage{nicefrac}       
\usepackage{microtype}      
\usepackage{amsmath}
\usepackage{amsthm}
\usepackage{graphicx}
\usepackage{algorithm}
\usepackage{algpseudocode}

\theoremstyle{definition}
\newtheorem{defn}{Definition}[section]

\newtheorem{prop}{Property}

\newcommand{\reals}{{\mbox{\bf R}}}

\newcommand{\symm}{{\mbox{\bf S}}}  

\newcommand{\Tr}{\mathop{\bf Tr}}
\newcommand{\diag}{\mathop{\bf diag}}

\newcommand{\Expect}{\mathop{\bf E{}}}

\newcommand{\vectorize}{\mathop{\bf vec}}

\newcommand{\eg}{{\it e.g.}}

\title{A Matrix Gaussian Distribution}

\author{
  Shane Barratt\\
  Department of Electrical Engineering\\
  Stanford University\\
  Stanford, CA 94305 \\
  \texttt{sbarratt@stanford.edu} \\
}

\begin{document}

\maketitle

\begin{abstract}
In this note, we define a Gaussian probability distribution over matrices.
We prove some useful properties of this distribution, namely, the fact that marginalization, conditioning, and affine transformations preserve the matrix Gaussian distribution.
We also derive useful results regarding the expected value of certain quadratic forms based solely on covariances between matrices.
Previous definitions of matrix normal distributions are severely under-parameterized, assuming unrealistic structure on the covariance (see Section 2).
We believe that our generalization is better equipped for use in practice.
\end{abstract}

\section{Introduction}

Recall that a random vector $x\in\reals^n$ has a Gaussian distribution if its probability distribution is fully characterized by a mean vector $\mu\in\reals^n$ and covariance matrix $\Sigma\in\symm_{++}^{n\times n}$.
If $x$ has a Gaussian distribution we write $x\sim\mathcal{N}(\mu,\Sigma)$.
The multivariate Gaussian distribution has many nice properties, \eg, its marginal and conditional distributions also have Gaussian distributions.
See Appendix A for a review of some of these properties.

A natural question to ask is: what is the best generalization of the multivariate Gaussian distribution to matrices?
To deal with matrix-valued distributions, we need to introduce two useful mathematical operators: the Kronecker product and vectorization. (Some of their useful properties are given in Appendix B.)

\begin{defn}[Kronecker product]
The Kronecker product between two matrices $A\in\reals^{m\times n}$ and $B\in\reals^{p\times q}$ is denoted $A\otimes B\in\reals^{mp\times nq}$ and is given by
\[A\otimes B = \begin{bmatrix}a_{11}B & \cdots & a_{1n}B \\
\vdots & \ddots & \vdots \\
a_{m1}B & \cdots & a_{mn}B\end{bmatrix}.\]
\end{defn}

\begin{defn}[Vectorization]
The vectorization of a matrix $A\in\reals^{m\times n}$ is denoted $\vectorize(A)\in\reals^{nm}$ and is equal to the concatenation of its columns, or
\[\vectorize(A)=\begin{bmatrix}a_1\\a_2\\ \vdots \\ a_n\end{bmatrix}\]
where $a_1,\ldots,a_n$ are the columns of $A$.
The inverse is given by $\vectorize{}^{-1}(A,m,n)$, where the second two arguments are needed to specify the size of the resulting matrix.
\end{defn}

Say we have a matrix $A\in\reals^{m\times n}$.
The most natural representation for $A$ as having a Gaussian distribution is to give $\vectorize(A)$ a (standard) multivariate Gaussian distribution and let $M\in\reals^{m\times n}$ and $\Sigma\in\reals^{nm\times nm}$ be the mean and covariance, or
\begin{equation*}
\begin{aligned}
M &= \Expect\left[A\right]\\
\Sigma &= \Expect\left[\vectorize(A-M)\vectorize(A-M)^T\right].
\end{aligned}
\end{equation*}

There is another representation based on the Kronecker product, where we let the mean $M\in\reals^{m\times n}$ and the covariance $S\in\reals^{m^2\times n^2}$, or
\begin{equation*}
\begin{aligned}
M &= \Expect\left[A\right]\\
S &= \Expect\left[(A-M)\otimes(A-M)\right]\text,
\end{aligned}
\end{equation*}
where the expectations are taken over $A$.
There is a mapping between the covariance in these two definitions, given by
\begin{equation*}
\begin{aligned}
\Sigma_{:,i+jm+1}&=\vectorize(S_{i(m):(i+1)m,j(n):(j+1)n})
\end{aligned}
\end{equation*}
where $\Sigma_{:,i}$ is the $i$th column of $\Sigma$.
Roughly, the vectorization of the $(i,j)$th block of $S$ corresponds to a column in $\Sigma$.
Because of this (albeit messy) one-to-one correspondence, we can swap between the two representations when it is convenient.
We will then write $A\sim\mathcal{N}(M,\Sigma)$ or $A\sim\mathcal{N}(M,S)$.
We now derive several useful properties of matrix-valued Gaussians.

\subsection{Probability Distribution and Entropy}
The distribution of $A\sim\mathcal{N}(M,\Sigma)$ is given by
\begin{equation*}
p(A) = \frac{1}{(2\pi)^{n/2}\det(\Sigma)^{1/2}}\exp\left\{-\frac{1}{2}\vectorize(A-M)^T\Sigma^{-1}\vectorize(A-M)\right\}.
\end{equation*}
Similarly, its differential entropy is given by $H(A)=\frac{1}{2}\ln\det(2\pi e\Sigma)$.
This is because of Property~\ref{prop:distr} in Appendix A.
Now, onto more interesting properties.

\subsection{Affine Transformations}
Let $A\in\reals^{m\times n}$ be such that $A\sim\mathcal{N}(M,\Sigma)$, and suppose $B\in\reals^{p\times m}$ and $C\in\reals^{m\times n}$ are constant matrices. Then the random matrix $BA+C$ follows a Gaussian distribution, or
\begin{equation*}
BA+C\sim\mathcal{N}(BM+C, (I\otimes B)\Sigma(I\otimes B)^T).
\end{equation*}
\begin{proof}
First, we note that
\[
\begin{aligned}
\vectorize(BA+C) &= \vectorize(BA)+\vectorize(C) \\
&= (I \otimes B) \vectorize(A) + \vectorize(C). \\
\end{aligned}
\]
where the first equality is because $\vectorize$ is linear (it is easy to convince yourself of this) and the second line is because of Property~\ref{prop:kron} in Appendix B.
By Property~\ref{prop:affine} in Appendix A, we know that $\vectorize(BA+C)$ follows a Gaussian distribution.
We then simply evaluate the mean
\[\Expect\left[BA+C\right] = BM+C\]
and the covariance
\[\text{Cov}[\vectorize(BA+C)]=(I \otimes B)\Sigma(I \otimes B)^T\text,\]
thus completing the proof.
\end{proof}

\subsection{Marginal and Conditional Distributions}

Let $[A \;B]\in\reals^{m\times (n_A + n_B)}$ and $[A\;B]\sim\mathcal{N}(\begin{bmatrix}M_A & M_B \end{bmatrix},\begin{bmatrix}\Sigma_{AA} & \Sigma_{AB} \\ \Sigma_{AB}^T & \Sigma_{BB}\end{bmatrix})$.

Then $A \sim \mathcal{N}(M_A,\Sigma_{AA})$ and $B\sim\mathcal{N}(M_B,\Sigma_{BB})$.
This follows directly from the marginal distribution properties of the multivariate Gaussian distribution.

What is the distribution of $A$ \emph{given} $B=C$ for some $C\in\reals^{m\times n_B}$?
Unsurprisingly, it is also a matrix Gaussian, given by
\begin{equation*}
\begin{aligned}
M_{A\mid B=C}&=M_A + \vectorize{}^{-1}(\Sigma_{AB}\Sigma_{BB}^{-1}\vectorize(C-M_B),m,n_A)\\
\Sigma_{A\mid B=C}&=\Sigma_{AA}-\Sigma_{AB}\Sigma_{BB}^{-1}\Sigma_{AB}^T.
\end{aligned}
\end{equation*}
This follows from Property~\ref{prop:conditional} in Appendix A.

\subsection{Expectations of Matrix Quadratic Forms}
Let $A\in\reals^{m\times n}$, $B\in\reals^{m\times p}$ are random matrices, and $C\in\reals^{m\times m}$ is a constant matrix.
Suppose we have the mean and covariance of $A$ and $B$, or
\begin{equation*}
\begin{aligned}
M_A &= \Expect\left[A\right]\\
M_B &= \Expect\left[B\right]\\
S_{BA} &= \Expect\left[(B-M_B)\otimes(A-M_A)\right].
\end{aligned}
\end{equation*}
Then
\[
\Expect[A^TCB]=\vectorize{}^{-1}(S_{BA}^T\vectorize(C),n,p)+M_A^TCM_B.
\]
This is a useful result.
Because expectation is linear, for fixed vectors $x\in\reals^n$ and $u\in\reals^p$, we can write the quadratic form $x^TA^TCBu$ as $x^T\Expect\left[A^TCB\right]u$.
The above result is quite general; for example, we can plug $B=A$ and get the result for $\Expect[A^TCA]$.
\begin{proof}
It is easy to see that
\[\Expect[A^TCB] = \Expect[(A-M_A)^TC(B-M_B)] + M_A^T C M_B\]
by using the fact that expectation is linear.
Then
\begin{equation*}
\begin{aligned}
\Expect[\vectorize((A-M_A)^TC(B-M_B))] &= \Expect[(B-M_B)^T\otimes(A-M_A)^T] \vectorize(C)\\
&= \Expect[(B-M_B)\otimes(A-M_A)]^T \vectorize(C)\\
&= S_{BA}^T \vectorize(C)
\end{aligned}
\end{equation*}
where the first equality is because of Property~\ref{prop:kron} in Appendix B and the second equality is because of Property~\ref{prop:kron-transpose} in Appendix B.
The result follows by performing the inverse of the vectorization operation.
\end{proof}

\section{Flaws with the ``Matrix Normal'' Distribution}

The ``matrix normal'' distribution was first introduced by De Waal in 1985~\cite{de1985matrix}.
Let $X\in\reals^{n\times p}$. Then $X$ follows the matrix normal distribution 
\[\mathcal{M}\mathcal{N}_{n,p}(M,U,V)\]
for $M\in\reals^{n\times p}$, $U\in\symm_{+}^{n\times n}$, and $V\in\symm_+^{p\times p}$ if and only if
\[\vectorize(X)\sim\mathcal{N}_{np}(\vectorize(M),V\otimes U).\]
The matrix normal distribution has several nice properties~\cite{gupta1999matrix}, but it is severely under-parameterized.
The matrix normal has only $n^2+p^2$ parameters for the covariance matrix, which technically has $(np)^2$ degrees of freedom.
Letting $n=p$, this means that the matrix normal distribution only has a fraction $2n^2/(n^2n^2)=2/n^2$ parameters.
For modest $n$, \eg, $n=10$, the matrix normal only has 2\% of the parameters.
This could be beneficial for storage purposes, however, by a simple two-dimensional counter-example, we will show that it fails to represent many common matrix Gaussian distributions that are easy for the matrix Gaussian distribution.

Suppose $n=p=2$ and $M=0$.
Also suppose that all the entries in $X$ are i.i.d. with the variance of $X_{ij}=\sigma_{ij}$, meaning that $V\otimes U$ should be equal to $\diag([\sigma_{11},\sigma_{12},\sigma_{21},\sigma_{22}])$.
This is a reasonable case for a prior distribution on $X$.
For this to be representable with the matrix normal distribution, we need to have that
\begin{equation*}
\begin{aligned}
V_{11}U_{11}&=\sigma_{11}\\
V_{11}U_{22}&=\sigma_{12}\\
V_{22}U_{11}&=\sigma_{21}\\
V_{22}U_{22}&=\sigma_{22}.\\
\end{aligned}
\end{equation*}
Dividing the 1st equation by the 2nd and the 3rd equation by the 4th, we get that
\begin{equation*}
\sigma_{11}/\sigma_{12} = \sigma_{21}/\sigma_{22}.
\end{equation*}
Therefore, the matrix normal in two dimensions can only represent i.i.d. covariances when the above equation holds for the priors, making it unusable in practice.




\subsubsection*{Acknowledgments}

This material is based upon work supported by the National Science Foundation Graduate Research Fellowship under Grant No. DGE-1656518.

\small

\bibliographystyle{unsrt}
\bibliography{mybib}

\appendix
\section*{Appendix A}
We review some properties of multivariate Gaussians. For proofs of these, see any statistics textbook, \eg,~\cite{murphy2013machine}.
\begin{prop}
\label{prop:distr}
Let $x\sim\mathcal{N}(\mu,\Sigma)$. Then the probability distribution of $x$ is given by
\[p(x) = \frac{1}{(2\pi)^{n/2}\det(\Sigma)^{1/2}}\exp\left\{-\frac{1}{2}(x-\mu)^T\Sigma^{-1}(x-\mu)\right\}\text,\]
and the differential entropy is given by
\[H(x) = \frac{1}{2}\ln\det(2\pi e\Sigma).\]
\end{prop}

\begin{prop}
\label{prop:affine}
Let $x\in\reals^n$ and $x\sim\mathcal{N}(\mu,\Sigma)$, $A\in\reals^{m\times n}$ and $b\in\reals^n$. Then 
$Ax+b\sim\mathcal{N}(A\mu,A\Sigma A^T)$.
\end{prop}

\begin{prop}
\label{prop:marginal}
Let $\begin{bmatrix}x \\ y\end{bmatrix}\sim\mathcal{N}(\begin{bmatrix}\mu_x \\ \mu_y \end{bmatrix},\begin{bmatrix}\Sigma_{xx} & \Sigma_{xy} \\ \Sigma_{xy}^T & \Sigma_{yy}\end{bmatrix})$. Then the marginal distributions of $x$ and $y$ are given by
\[x\sim\mathcal{N}(\mu_x, \Sigma_{xx}).\]
and
\[y\sim\mathcal{N}(\mu_y, \Sigma_{yy}).\]
\end{prop}

\begin{prop}
\label{prop:conditional}
Let $\begin{bmatrix}x \\ y\end{bmatrix}\sim\mathcal{N}(\begin{bmatrix}\mu_x \\ \mu_y \end{bmatrix},\begin{bmatrix}\Sigma_{xx} & \Sigma_{xy} \\ \Sigma_{xy}^T & \Sigma_{yy}\end{bmatrix})$. Then the conditional distribution of $x$ \emph{given} $y=a$ is given by $\mathcal{N}(\mu_{x\mid y=a},\Sigma_{x\mid y=a})$ where
\[\mu_{x\mid y=a}=\mu_x + \Sigma_{xy}\Sigma_{yy}^{-1}(a-\mu_y)\]
and
\[\Sigma_{x\mid y=a}=\Sigma_{xx}-\Sigma_{xy}\Sigma_{yy}^{-1}\Sigma_{xy}^T.\]
\end{prop}

\begin{prop}
\label{prop:quadform}
Let $x\in\reals^n$ and $x\sim\mathcal{N}(\mu,\Sigma)$ and $A\in\reals^{n\times n}$. Then the expectation of the quadratic form $x^TAx$ is given by
\[\Expect\left[x^TAx\right] = \Tr(A\Sigma)+\mu^TA\mu.\]
\end{prop}

\section*{Appendix B}
We review some properties of Kronecker products, vectorizations, and traces.
For a proof of these, see~\cite{graham1981kronecker}.

\begin{prop}
\label{prop:kron-transpose}
Given $A\in\reals^{m\times n}$ and $B\in\reals^{n\times p}$, we have that
\[(A^T \otimes B^T) = (A\otimes B)^T.\]
\end{prop}

\begin{prop}
\label{prop:kron}
Given $A\in\reals^{m\times n}$, $B\in\reals^{n\times p}$, and $C\in\reals^{p\times k}$, we have that
\[\vectorize(ABC)=(C^T\otimes A)\vectorize(B).\]
\end{prop}

\begin{prop}
\label{prop:tr}
Given $A\in\reals^{m\times n}$, $B\in\reals^{n\times p}$, and $C\in\reals^{p\times k}$, we have that
\[\Tr(A^TB)=\vectorize(A)^T\vectorize(B).\]
\end{prop}

\end{document}